%Paper: funct-an/9411002
%From: Graziano Crasta - SISSA Trieste <crasta@tsmi19.sissa.it>
%Date: Thu, 03 Nov 1994 12:03:47 +0200
%Date (revised): Wed, 16 Nov 1994 15:53:25 +0200

%
% CRAMAL
%
% 16/09/94
%
% ========================
%
\magnification=\magstep1
\def\autori{G.~Crasta and A.~Malusa}
\def\titolo{Existence results for non--coercive variational problems}
\newif\iftitlepage
\newif\ifproofmode
\proofmodefalse
\titlepagetrue
%
%\input macro.tex
%
% MACRO.TEX
% ========================
%
% FONTS
%
\font\twelverm=cmr12
\font\twelvei=cmmi12
\font\twelvesy=cmsy10
\font\twelvebf=cmbx12
\font\twelvett=cmtt12
\font\twelveit=cmti12
\font\twelvesl=cmsl12
\font\ninerm=cmr9
\font\ninei=cmmi9
\font\ninesy=cmsy9
\font\ninebf=cmbx9
\font\ninett=cmtt9
\font\nineit=cmti9
\font\ninesl=cmsl9
\font\ninesc = cmcsc10 at 9pt

\font\eightrm=cmr8
\font\eighti=cmmi8
\font\eightsy=cmsy8
\font\eightbf=cmbx8
\font\eighttt=cmtt8
\font\eightit=cmti8
\font\eightsl=cmsl8
\font\sixrm=cmr6
\font\sixi=cmmi6
\font\sixsy=cmsy6
\font\sixbf=cmbx6
\catcode`@=11 % we will access private macros of plain TeX (carefully)
\newskip\ttglue
%
%MACRO TWELVEPOINT
\def\twelvepoint{\def\rm{\fam0\twelverm}% switch to 12-point type
\textfont0=\twelverm  \scriptfont0=\ninerm
\scriptscriptfont0=\sevenrm
\textfont1=\twelvei  \scriptfont1=\ninei  \scriptscriptfont1=\seveni
\textfont2=\twelvesy  \scriptfont2=\ninesy
\scriptscriptfont2=\sevensy
\textfont3=\tenex  \scriptfont3=\tenex  \scriptscriptfont3=\tenex
\textfont\itfam=\twelveit  \def\it{\fam\itfam\twelveit}%
\textfont\slfam=\twelvesl  \def\sl{\fam\slfam\twelvesl}%
\textfont\ttfam=\twelvett  \def\tt{\fam\ttfam\twelvett}%
\textfont\bffam=\twelvebf  \scriptfont\bffam=\ninebf
\scriptscriptfont\bffam=\sevenbf  \def\bf{\fam\bffam\twelvebf}%
\tt  \ttglue=.5em plus.25em minus.15em
\normalbaselineskip=15pt
\setbox\strutbox=\hbox{\vrule height10pt depth5pt width0pt}%
\let\sc=\tenrm  \let\big=\twelvebig  \normalbaselines\rm}
%
%MACRO TENPOINT
\def\tenpoint{\def\rm{\fam0\tenrm}% switch to 10-point type
\textfont0=\tenrm  \scriptfont0=\sevenrm  \scriptscriptfont0=\fiverm
\textfont1=\teni  \scriptfont1=\seveni  \scriptscriptfont1=\fivei
\textfont2=\tensy  \scriptfont2=\sevensy  \scriptscriptfont2=\fivesy
\textfont3=\tenex  \scriptfont3=\tenex  \scriptscriptfont3=\tenex
\textfont\itfam=\tenit  \def\it{\fam\itfam\tenit}%
\textfont\slfam=\tensl  \def\sl{\fam\slfam\tensl}%
\textfont\ttfam=\tentt  \def\tt{\fam\ttfam\tentt}%
\textfont\bffam=\tenbf  \scriptfont\bffam=\sevenbf
\scriptscriptfont\bffam=\fivebf  \def\bf{\fam\bffam\tenbf}%
\tt  \ttglue=.5em plus.25em minus.15em
\normalbaselineskip=12pt
\setbox\strutbox=\hbox{\vrule height8.5pt depth3.5pt width0pt}%
\let\sc=\eightrm  \let\big=\tenbig  \normalbaselines\rm}
%
%MACRO NINEPOINT
\def\ninepoint{\def\rm{\fam0\ninerm}% switch to 9-point type
\textfont0=\ninerm  \scriptfont0=\sixrm  \scriptscriptfont0=\fiverm
\textfont1=\ninei  \scriptfont1=\sixi  \scriptscriptfont1=\fivei
\textfont2=\ninesy  \scriptfont2=\sixsy  \scriptscriptfont2=\fivesy
\textfont3=\tenex  \scriptfont3=\tenex  \scriptscriptfont3=\tenex
\textfont\itfam=\nineit  \def\it{\fam\itfam\nineit}%
\textfont\slfam=\ninesl  \def\sl{\fam\slfam\ninesl}%
\textfont\ttfam=\ninett  \def\tt{\fam\ttfam\ninett}%
\textfont\bffam=\ninebf  \scriptfont\bffam=\sixbf
\scriptscriptfont\bffam=\fivebf  \def\bf{\fam\bffam\ninebf}%
\tt  \ttglue=.5em plus.25em minus.15em
\normalbaselineskip=11pt
\setbox\strutbox=\hbox{\vrule height8pt depth3pt width0pt}%
\let\sc=\sevenrm  \let\big=\ninebig  \normalbaselines\rm}
%
%MACRO EIGHTPOINT
\def\eightpoint{\def\rm{\fam0\eightrm}% switch to 8-point type
\textfont0=\eightrm  \scriptfont0=\sixrm  \scriptscriptfont0=\fiverm
\textfont1=\eighti  \scriptfont1=\sixi  \scriptscriptfont1=\fivei
\textfont2=\eightsy  \scriptfont2=\sixsy  \scriptscriptfont2=\fivesy
\textfont3=\tenex  \scriptfont3=\tenex  \scriptscriptfont3=\tenex
\textfont\itfam=\eightit  \def\it{\fam\itfam\eightit}%
\textfont\slfam=\eightsl  \def\sl{\fam\slfam\eightsl}%
\textfont\ttfam=\eighttt  \def\tt{\fam\ttfam\eighttt}%
\textfont\bffam=\eightbf  \scriptfont\bffam=\sixbf
\scriptscriptfont\bffam=\fivebf  \def\bf{\fam\bffam\eightbf}%
\tt  \ttglue=.5em plus.25em minus.15em
\normalbaselineskip=9pt
\setbox\strutbox=\hbox{\vrule height7pt depth2pt width0pt}%
\let\sc=\sixrm  \let\big=\eightbig  \normalbaselines\rm}
%
%MACRO BIG
\def\twelvebig#1{{\hbox{$\textfont0=\twelverm\textfont2=\twelvesy
	\left#1\vbox to10pt{}\right.\n@space$}}}
\def\tenbig#1{{\hbox{$\left#1\vbox to8.5pt{}\right.\n@space$}}}
\def\ninebig#1{{\hbox{$\textfont0=\tenrm\textfont2=\tensy
	\left#1\vbox to7.25pt{}\right.\n@space$}}}
\def\eightbig#1{{\hbox{$\textfont0=\ninerm\textfont2=\ninesy
	\left#1\vbox to6.5pt{}\right.\n@space$}}}

%for 10-point math in 9-point territory
\font\medbf=cmbx10 scaled\magstep2
\def\today{\ifcase\month\or
January\or February\or March\or April\or May\or June\or
July\or August\or September\or October\or November\or December\fi
\space\number\day, \number\year}
%
%
% ================================================================
%
% Definizione delle macro
%
\nopagenumbers
\def\proof{\noindent{\bf Proof.}\quad}
\def\sqr#1#2{\vbox{
   \hrule height .#2pt
   \hbox{\vrule width .#2pt height #1pt \kern #1pt
      \vrule width .#2pt}
   \hrule height .#2pt}}

\def\finedim{\hfill\hbox{$\sqr74$} \goodbreak\medskip}

%
% MACROS FOR BIBLIOGRAPHY
%
% \bibart{label}{author}{title}{review}{vol}{year}{start-page}{end-page}
% \bibprep{label}{author}{title}{type}{year}
% \biblib{label}{author}{title}{editor}{place}{year}
% \bibinlib{label}{author}{art-title}{book-title}{editor}{place}{year}
% \bibproc{label}{author}{title}{proceedings}{place}{year}
% \bibtoapp{label}{author}{title}
%
\newcount\vol \newcount\pag
\def\bibart#1#2#3#4#5#6#7#8{\global\vol=#5 \global\pag=#7
%\par\hangafter=10
{\item{[{\bf\bib{#1}}]}{\ninesc #2}: {#3},
{\sl #4},%
{\ifnum\vol=0\else{{\bf\ #5},}\fi}
{#6}%
{\ifnum\pag=0\else{, p.~{#7}--{#8}}\fi}%
. \medskip}}
\def\bibprep#1#2#3#4#5{
%\par\hangafter=10
{\item{[{\bf\bib{#1}}]}{\ninesc #2}: {#3},
{\sl #4}, {#5}. \medskip}}
\def\biblib#1#2#3#4#5#6{
%\par\hangafter=10%
{\item{[{\bf\bib{#1}}]}{\ninesc #2}: {\sl #3}, {#4}, {#5}, {#6}.
\medskip}}
\mathchardef\emptyset="001F
%
%   MACROS FOR LABELLING
%
%   In what follows, <label> or <plabel> should be replaces by actual labels.
%   \chp        \head \chp{<label>}. Preliminaries \endhead
%   \thm        \proclaim{\thm{<label>}. Theorem}
%   \frm 	\tag{\frm{<label>}}
%   \bib        \ref\no\bib{<label>}   \by .. \paper .. \jour ...\endref
%   \rf         By  (\rf{<label>}) we get...
%               or, more ortodoxly,
%               By  \rom{(\rf{<label>})} we get...
%   \pgp        This is an important place \pgp{<plabel>}.
%   \rfp        This can be found on p. \rfp{<plabel>}.
%   \ix2        We define \ix2{<plabel>}{gradient} {\it gradient\/} by ...
%   \ix3        We introduce the space \ix3{<plabel>}{$W^{1,p}$}
%               $W^{1,p}(\Omega}$.
%   \newpar     \newpar{label}{title of paragraph}
%   \nfor       \nfor{label}
%   \theo
%   \lemma
%   \prop
%   \corol
%   \dhef       \xxx{label}{statement}
%   \ehse       \ehse{label}
%   \ohss       \ohss{label}
%
\font\sixrm=cmr6
\newcount\tagno \tagno=0		        % equation #'s
\newcount\thmno	\thmno=0	         	% theorem numbers
\newcount\bibno	\bibno=0			% bibliography numbers
\newcount\chapno\chapno=0                       % chapters
\newcount\verno            % number of the version
\newif\ifwanted
\wantedfalse
\newif\ifindexed
\indexedfalse
% AUXILIARY MACROS
\def\ifundefined#1{\expandafter\ifx\csname+#1\endcsname\relax}
\def\Wanted#1{\ifundefined{#1} \wantedtrue \immediate\write0{Wanted #1
\the\chapno.\the\thmno}\fi}
\def\Increase#1{{\global\advance#1 by 1}}
%for our numbering
\def\Assign#1#2{\immediate
\write1{\noexpand\expandafter\noexpand\def
 \noexpand\csname+#1\endcsname{#2}}\relax
 \global\expandafter\edef\csname+#1\endcsname{#2}}
%for pages
\def\pAssign#1#2{\write1{\noexpand\expandafter\noexpand\def
 \noexpand\csname+#1\endcsname{#2}}}
\def\lPut#1{\ifproofmode\llap{\hbox{\sixrm #1\ \ \ }}\fi}
\def\rPut#1{\ifproofmode$^{\hbox{\sixrm #1}}$\fi}
%
% USER'S MACROS
% Number of a new chapter labelled as #1:
\def\chp#1{\global\tagno=0\global\thmno=0\Increase\chapno
\Assign{#1}
{\the\chapno}{\lPut{#1}\the\chapno}}
% Paragraph: \newpar{label}{title}
\def\newpar#1#2{\bigskip {\noindent\bf \chp{#1}. #2} \medskip}  %
% Number of a  new theorem (remark, definition, example,...) labelled as #1:
\def\thm#1{\Increase\thmno
\Assign{#1}{\the\chapno.\the\thmno}\the\chapno.\the\thmno\rPut{#1}}
% Number of a new formula labelled as #1:
\def\frm#1{\Increase\tagno
  \Assign{#1}{\the\chapno.\the\tagno}\lPut{#1}{\the\chapno.\the\tagno}}
\def\nfor#1{\leqno{(\frm{#1})}}
% Number of a new bibliographical item labelled as #1:
\def\bib#1{\Increase\bibno
\Assign{#1}{\the\bibno}\lPut{#1}{\the\bibno}}
% Remembers the page number of the place labelled as #1
\def\pgp#1{\pAssign{#1/}{\the\pageno}}
% Creates the index item #3 labelled as #2: it is \ix2 for
% index and \ix3 for  list of symbols.
\def\ix#1#2#3{\pAssign{#2}{\the\pageno}
\immediate\write#1{\noexpand\idxitem{#3}{\noexpand\csname+#2\endcsname}}}
% Number of the item labelled as #1, for references:
\def\rf#1{\Wanted{#1}\csname+#1\endcsname\relax\rPut {#1}}
% Page number of the place labelled as #1, for references:
\def\rfp#1{\Wanted{#1}\csname+#1/\endcsname\relax\rPut{#1}}
% Theorems, propositions, etc.
\long\def\theo#1#2{\medskip {\noindent\bf Theorem \thm{#1}}
\quad{\sl #2}\medskip}
\long\def\prop#1#2{\medskip {\noindent\bf
Proposition \thm{#1}}\quad{\sl #2}\medskip}
\long\def\lemma#1#2{\medskip {\noindent\bf
Lemma \thm{#1}}\quad{\sl #2}\medskip}
\long\def\corol#1#2{\medskip {\noindent\bf
Corollary \thm{#1}}\quad{\sl #2}\medskip}
\long\def\dhef#1#2{\medskip {\noindent\bf
Definition \thm{#1}}\quad{\sl #2}\medskip}

\def\ohss#1{\medskip {\noindent\bf Remark \thm{#1}}\quad}
%
%\newread\fileaux
%\openin\fileaux=\jobname.aux
%\ifeof\fileaux
%\verno=0
%\else
%\closein\fileaux
%\input \jobname.aux
%\fi
%\Increase\verno
\verno =7
\expandafter \def \csname +intro\endcsname {1}
\expandafter \def \csname +mina\endcsname {1.1}
\expandafter \def \csname +enex\endcsname {1.2}
\expandafter \def \csname +mind\endcsname {1.3}
\expandafter \def \csname +prelim\endcsname {2}
\expandafter \def \csname +coA\endcsname {2.1}
\expandafter \def \csname +cosb\endcsname {2.2}
\expandafter \def \csname +subgr\endcsname {2.3}
\expandafter \def \csname +propsb\endcsname {2.1}
\expandafter \def \csname +clo\endcsname {3}
\expandafter \def \csname +SCI\endcsname {3.1}
\expandafter \def \csname +caratt\endcsname {3.2}
\expandafter \def \csname +claF\endcsname {3.3}
\expandafter \def \csname +remuno\endcsname {3.4}
\expandafter \def \csname +luno\endcsname {3.5}
\expandafter \def \csname +clos\endcsname {3.6}
\expandafter \def \csname +fca\endcsname {3.1}
\expandafter \def \csname +fcaa\endcsname {3.2}
\expandafter \def \csname +fcb\endcsname {3.3}
\expandafter \def \csname +aab\endcsname {3.4}
\expandafter \def \csname +fcc\endcsname {3.5}
\expandafter \def \csname +fcd\endcsname {3.6}
\expandafter \def \csname +aac\endcsname {3.7}
\expandafter \def \csname +coruno\endcsname {3.7}
\expandafter \def \csname +cordue\endcsname {3.8}
\expandafter \def \csname +exi\endcsname {4}
\expandafter \def \csname +claE\endcsname {4.1}
\expandafter \def \csname +cen\endcsname {4.1}
\expandafter \def \csname +carE\endcsname {4.2}
\expandafter \def \csname +cend\endcsname {4.2}
\expandafter \def \csname +sta\endcsname {4.3}
\expandafter \def \csname +stb\endcsname {4.4}
\expandafter \def \csname +stc\endcsname {4.5}
\expandafter \def \csname +std\endcsname {4.6}
\expandafter \def \csname +remzero\endcsname {4.3}
\expandafter \def \csname +EG\endcsname {4.4}
\expandafter \def \csname +esti\endcsname {4.5}
\expandafter \def \csname +hfitre\endcsname {4.7}
\expandafter \def \csname +beuno\endcsname {4.8}
\expandafter \def \csname +bedue\endcsname {4.9}
\expandafter \def \csname +nag\endcsname {4.6}
\expandafter \def \csname +bound\endcsname {4.7}
\expandafter \def \csname +app\endcsname {4.8}
\expandafter \def \csname +min\endcsname {4.10}
\expandafter \def \csname +DR\endcsname {4.11}
\expandafter \def \csname +pteta\endcsname {4.12}
\expandafter \def \csname +unb\endcsname {4.13}
\expandafter \def \csname +DRk\endcsname {4.14}
\expandafter \def \csname +bdaa\endcsname {4.15}
\expandafter \def \csname +cip\endcsname {4.16}
\expandafter \def \csname +miao\endcsname {4.17}
\expandafter \def \csname +posi\endcsname {4.18}
\expandafter \def \csname +bau\endcsname {4.19}
\expandafter \def \csname +enin\endcsname {4.20}
\expandafter \def \csname +ben\endcsname {4.21}
\expandafter \def \csname +cca\endcsname {4.9}
\expandafter \def \csname +ccfa\endcsname {4.22}
\expandafter \def \csname +ccfb\endcsname {4.23}
\expandafter \def \csname +ccfc\endcsname {4.24}
\expandafter \def \csname +ccr\endcsname {4.10}
\expandafter \def \csname +ccb\endcsname {4.11}
\expandafter \def \csname +conca\endcsname {4.12}
\expandafter \def \csname +AAB\endcsname {1}
\expandafter \def \csname +BD\endcsname {2}
\expandafter \def \csname +BM\endcsname {3}
\expandafter \def \csname +Butt\endcsname {4}
\expandafter \def \csname +CC\endcsname {5}
\expandafter \def \csname +CTZ\endcsname {6}
\expandafter \def \csname +ces\endcsname {7}
\expandafter \def \csname +Cl\endcsname {8}
\expandafter \def \csname +Clar\endcsname {9}
\expandafter \def \csname +Cla\endcsname {10}
\expandafter \def \csname +CL\endcsname {11}
\expandafter \def \csname +Cra\endcsname {12}
\expandafter \def \csname +ET\endcsname {13}
\expandafter \def \csname +Mar\endcsname {14}
\expandafter \def \csname +Ol\endcsname {15}
\expandafter \def \csname +Ray\endcsname {16}

\immediate\openout1=\jobname.aux
\immediate\write1{\noexpand\verno=\the\verno}
\ifindexed
\immediate\openout2=\jobname.idx
\immediate\openout3=\jobname.sym
\fi
%To be set in front of the end of the job: \prebye\enddocument
\def\prebye{\ifwanted
\message{Warning: Undefined references! Rerunning could help}\fi}
\headline={\ifnum\pageno>0\ifodd\pageno\rightheadline
\else\leftheadline\fi\fi}
\def\rightheadline{\hfil{\eightpoint\titolo}
\hfil\tenrm\folio}
\def\leftheadline{\tenrm\folio\hfil{\eightpoint\autori}
\hfil} \topskip=25pt
%
% =================================================================
%
% Macros for typesetting the text
%
\def \nat{I\!\!N}

\def\R{{I\!\!R}}

\def \C{{\cal C}}
\def\D{{\cal D}}

\def\E{{\cal E}}
\def\G{{\cal G}}

\def\dug{{\ \buildrel .\over =\ }}

\let\meti=\item
%
%
% Altre definizioni
%
\def\coepif{{\rm co\,epi\,f}}
\def\coepifi{{\rm co\,epi\,\phi}}

\def\D{{\cal D}}

\def\E{{\cal E}}
\def\G{{\cal G}}
\def\fs{f^{**}}
\def\fis{\phi^{**}}
\def\sumj{\sum_{j=1}^{m+2}}
\def\sumI{\sum_{j\in I}}
\def\sumJ{\sum_{j\in J}}
\def\xianj{(\xi^n_j, a^n_j)}
\def\oxij{\overline{\xi}_j}
\def\olaj{{\lambda_j}}
\def\suth{\ \big|\ }
\def\limn{\lim_{n\rightarrow+\infty}}
\def\vmu{s}
\def\vfx{\varphi_{\nu,\xi}}
\def\psx{\psi_{\nu,\xi}}

\def\ptf{\partial_t\varphi}
\def\scal#1#2{\langle #1\,,\,#2\rangle}
\def\Rex{(-\infty,+\infty]}
\def\Rbar{\overline{\R}}
\def\dom{{\rm dom}}
\def\intdom#1{{\rm int(dom(}#1{\rm))}}
\def\dexi{\partial_{\xi}}
\def\psin{\psi^{**}(t^n,\xi^n)}
\def\psinj{\psi^{**}(t^n,\xi^n_j)}
\def\sj{\sum_{j=1}^{m+1}}
\def\lnj{\lambda^n_j}
\def\ect{e^{C_0 T}}
%
%=============================================
%
% TITLE PAGE
%
\def\fs{f^{**}}
\def\titlea{Existence results for non--coercive}
\def\titleb{variational problems}
\def\abstract{
We consider the problem (P):
$$\min\left\{\int_0^T\left [g(t,u)+f(t,u')\right] dt
\suth u\in W^{1,1}
([0,T], \R^m),\, u(0)=a,\, u(T)=b\right\}
$$
with neither coercivity nor convexity assumptions.
More precisely, assuming that $f(t,\xi)$
is a lower semicontinuous function, bounded from below,
Lipschitz continuous with respect to $t$,
satisfying
$$\limn\left[\fs(t^n,\xi^n)-
\scal{\nabla\fs(t^n,\xi^n)}{\xi^n}\right]=-\infty$$
for every sequence $\{t^n\}\in [0,T]$
and for every choice of points $\xi^n$ of differentiability
of $\fs(t^n,\cdot)$
such that $\lim_n|\xi^n|=+\infty$,
and assuming that $g(t,x)$ is a continuous function,
Lipschitz continuous with respect to $t$ and
concave with respect to $x$, such that
$g(t,x)\geq -\alpha -\beta |x|$
for every $(t,x)\in [0,T]\times \R^m$,
and for suitable constants $\alpha$ and $\beta\geq 0$,
we show that the problem~(P) has a solution in the space
$W^{1,\infty}([0,T],\R^m)$.
The main tools for the proof
are an existence theorem
for the problem~(P), with $f(t,\xi)$ convex with respect to $\xi$
and $g$ continuous,
and the closure of the convex hull
of the epigraph of the function $f(t,\cdot)$.
}
\def\keywords{Existence theory, non--convex problems,
              non--coercive problems.}
\def\amsclass{49J05.}
\null
\pageno=0
\tenpoint
\baselineskip=15pt
{  % start shrink
\leftskip=2truecm
\rightskip=2truecm
\vskip 2 truecm
{\medbf
   \centerline{\titlea}
   \vskip 4 truemm
   \centerline{\titleb}
}
\vskip 8 truemm
\iftitlepage
   \centerline{Graziano CRASTA$^{(1)(2)}$ and
   Annalisa MALUSA$^{(1)(3)}$}
   \vskip 22 truemm
\else
   \centerline{\today}
   \vskip 10 truemm
\fi
\iftitlepage
 {\eightpoint % start eightpoint
   \noindent{\bf Abstract. }
   \abstract
   \medskip
   \noindent{\bf AMS Subject Classification. }
   \amsclass
   %\vskip 15truemm
   \vfill
   \noindent
   \item{(1)} S.I.S.S.A., Via Beirut 2, 34014 Trieste (Italy).

   \item{(2)} Dipartimento di Matematica Pura ed Applicata,
   Via Campi 213/b, 41100 Modena (Italy).

   \item{(3)} Istituto di Matematica, Facolt\`a di Architettura,
   via Monteoliveto 3, 80134 Napoli (Italy).

 } % end ninepoint
 %\noindent{\bf Key Words.}
 %\keywords
 %\vskip 2.5 truecm
 %\centerline{Mathematics Subject Classification: 34H05, 49A36}
 \vskip 1 truecm
 \centerline{Ref. S.I.S.S.A. 130/94/M, September 1994.}
\fi
}  % end shrink
\iftitlepage
%\vfill
\eject
\else
\vskip 5truemm
\fi

\headline={\ifnum\pageno>0\ifodd\pageno\rightheadline
\else\leftheadline\fi\fi}
\def\rightheadline{\hfil{\eightpoint\titolo}
\hfil\tenrm\folio}
\def\leftheadline{\tenrm\folio\hfil{\eightpoint\autori}
\hfil} \topskip=25pt
%
% ================================================================
%
% TESTO
%
%\input cmtesto.tex
%
% CMTESTO
%
% 30/10/94
%
% ========================
%
%%%%%%%
\newpar{intro}{Introduction}

It is well known that, if
$L$ is a continuous function, such that
$\xi\mapsto L(t,x,\xi)$ is convex and superlinear, then
the variational problem
$$\min\left\{\int_0^T L(t,u,u') dt\suth
u\in W^{1,1}([0,T],\R^m),\, u(0)=a,\, u(T)=b\right\},\nfor{mina}$$
has a solution (see for instance [\rf{ces}]).

In recent years, the possibility of avoiding the convexity
or the superlinearity assumption
was investigated by many authors.

Some existence results for non--convex coercive problems
were obtained in the case
$L(t,x,\xi)=g(t,x)+f(t,\xi)$
(see for instance [\rf{CC}], [\rf{Mar}], [\rf{Ray}] and
the references therein).
In particular, in [\rf{CC}] it was proved that the convexity
assumption on $f(t,\cdot)$ can be replaced
by the condition of concavity of $g(t,\cdot)$.

More recently, some techniques were developed in order to
treat convex but non--coercive problems.
In this case,
even if the functionals considered are lower semicontinuous
in the weak topology of $W^{1,1}([0,T],\R^m)$,
the direct method of the Calculus of Variations
can not be applied, due to the lack of compactness of the
minimizing sequences.

In [\rf{Cla}], it was studied the problem~(\rf{mina})
with $L$ continuous, bounded from below and
convex with respect to $\xi$, the superlinearity
being replaced by a weaker condition
which permits to construct a
relatively compact minimizing sequence,
obtained by considering
the minima of suitable coercive approximating problems.
The main step in the proof of the existence result
in [\rf{Cla}] was to show that every minimum point
of the approximating problems solves a generalized
DuBois--Reymond condition, which implies that
the minimizing sequence is bounded in
the space $W^{1,\infty}([0,T],\R^m)$.

A similar approach was used in [\rf{CTZ}]
for the autonomous problem  $L(t,x,\xi)=g(x)+f(\xi)$,
where $g$ is a nonnegative continuous function,
and $f\in\C^1(\R^m,\R)$ is a  strictly convex function
bounded from below, such that
$$\lim_{|\xi|\rightarrow+\infty}
\left[f(\xi)-\scal{\nabla f(\xi)}{\xi}\right]=
-\infty.\nfor{enex}$$
In that paper, it was proved that, for every rectifiable
curve $C$ in $\R^m$ joining $a$ to $b$ there exists
a unique solution to the problem~(\rf{mina})
restricted to the class
of all absolutely continuous parameterizations
$u\colon I\rightarrow\R^m$ of $C$.
Thus, every element $u_n$ of a minimizing sequence
can be replaced
by the minimum corresponding
to the curve parameterized by $u_n$.
It can be shown, still using a DuBois--Reymond condition
satisfied by those minima, and by (\rf{enex}),
that this new sequence is
bounded in $W^{1,\infty}([0,T],\R^m)$,
so that there exists a minimum point for (\rf{mina})
in this space.

In [\rf{Cra}] both the superlinearity and the convexity
assumptions were dropped
for lagrangians of the form
$L(t,x,\xi)=\scal{a(t)}{x}+f(\xi)$
where $f$ is a lower semicontinuous function
whose convexification $\fs$ satisfies (\rf{enex})
for every diverging sequence of points of differentiability of
the Lipschitz continuous function $\fs$.
The existence of a minimum is proved by a technique
relying only
on a Lyapunov type theorem due to Olech (see [\rf{Ol}]).

For other results concerning non--coercive problems
we mention [\rf{AAB}], [\rf{BD}] and [\rf{BM}].

In this paper we consider non--autonomous problems of the form
$$
\min\left\{\int_0^T\left [g(t,u)+f(t,u^\prime)\right]\, dt\suth
u\in W^{1,1}
([0,T], \R^m)\,,\ u(0)=a\,,\ u(T)=b\right\} \nfor{mind}
$$
with neither coercivity nor convexity assumptions.
More precisely, we introduce the class $\E$ of all functions
$\psi:[0,T]\times\R^m\rightarrow\R$, bounded from below,
such that $\psi(\cdot, \xi)$ is Lipschitz continuous for
every fixed $\xi\in\R^m$, $\psi(t,\cdot)$ is
lower semicontinuous and satisfies
$$\limn\left[\psin-
\scal{\nabla\psin}{\xi^n}\right]=-\infty$$
for every sequence $\{t^n\}\in [0,T]$
and for every choice of points $\xi^n$ of differentiability
of $\psi^{**}(t^n,\cdot)$
such that $\lim_n|\xi^n|=+\infty$.
We show that, if $f\in\E$ and there exists two constants
$A$ and $B$, $B>0$ such that $f(t,\xi)\ge -A+B|\xi|$
for every $(t,\xi)\in [0,T]\times\R^m$,
and $g(t,x)$ is a continuous function,
Lipschitz continuous with respect to $t$,
concave with respect to $x$, satisfying
$g(t,x)\geq -\alpha -\beta |x|$
for every $(t,x)\in[0,T]\times\R^m$, and
for suitable constants $\alpha$ and $0\le\beta \le B/T$,
then the problem~(\rf{mind}) has a solution in the space
$W^{1,\infty}([0,T],\R^m)$.
This result is the analogue for a class of non--coercive
functionals of the one in [\rf{CC}], but it is not a generalization
of that result, due to the additional requirement of
the Lipschitz continuity of the lagrangian with respect to
the variable $t$.
However this extra regularity allows us to obtain the
necessary conditions that, used at an intermediate step,
give also a regularity result, interesting by itself.

As a first step we prove an existence result for (\rf{mind})
requiring that $f$ be convex with respect to $\xi$ and
dropping the concavity assumption on $g$.
This can be done
following [\rf{Cla}] and making
suitable changes, due to the
the fact that the lagrangian is not bounded from below.
The second step, linking the convex to the non--convex case,
is based on a result concerning the
closure of the convex hull of the epigraph of functions
whose convexification is strictly convex at infinity
(that is, the graph of the convexification contains no rays).
This result is an extension of the classical theorem that holds
for superlinear functions (see [\rf{ET}]).
We want to remark that the notion of strict convexity at
infinity was still used in [\rf{CL}] in order to study non--coercive
problems of the type~(\rf{mina}) with the additional state
constraint $\|u\|_{L^\infty}<R$.
We shall prove that every function in the class $\E$
is strictly convex at infinity for every fixed $t$,
so that, by using the previous results and
the Lyapunov theorem on the range of non--atomic measures,
the existence result for the non--convex problems follows.
The regularity of the solution of~(\rf{mind}) is a consequence
of the regularity of the solution to the relaxed problem.

\smallskip\noindent
{\bf Acknowledgements.} The authors wish to thank
Professor Arrigo Cellina for kindly suggesting the problem.

\newpar{prelim}{Preliminaries}

%If $x\in\R^m$, we denote by
%$B_r(x)$ the open ball centered
%at $x$ with radius $r$, and with
%$\overline{B}_r(x)$ its closure.
%We shall use the notation $B_r$ for $B_r(0)$.
%We shall denote by $2^{\R^m}$ the family of all subset of $\R^m$, and
%$2^{\R^m}\backslash\emptyset$
%will mean $2^{\R^m}\backslash\{\emptyset\}$.

%A set $A\subset\R^m$ is convex if, for every $x,y\in A$ and
%for every $\lambda\in[0,1]$, we have:
%$$\lambda x+ (1-\lambda) y\in A.$$

We shall denote by $\scal xy$ the
standard scalar product of two vectors $x$, $y\in\R^m$.
For every $1\leq p\leq+\infty$,
we shall denote by $L^p(I,\R^m)$ and $W^{1,p}(I,\R^m)$,
respectively, the usual Lebesgue and Sobolev spaces
of functions from the interval $I\dug [0,T]$ into $\R^m$.
We shall use the symbol $\|\cdot\|_{L^p}$ to denote the
norm in $L^p(I,\R^m)$.

If $A\subset\R^m$, we shall denote by
${\rm int}\, A$ the interior of $A$, and by
${\rm co}\, A$ the convex hull of a $A$,
that is, the smallest convex set which
contains $A$.
It is well known that,
by Carath\'eodory's theorem,
the convex hull of $A$ can be characterized by
$${\rm co}\,A = \left\{x\in\R^m\suth  x=\sum_{i=1}^{m+1}\lambda_i x_i\,,\
\tilde{\lambda}\in E_{m+1}\,,\ x_i\in A,\ i=1,\ldots,m+1\right\},
\nfor{coA}$$
where $\tilde{\lambda}\dug (\lambda_1,\ldots\lambda_{m+1})$, and
$E_{m+1}$ denotes the standard simplex:
$$E_{m+1}\dug \left\{(\lambda_1,\ldots\lambda_{m+1})\in\R^{m+1}\suth
\lambda_i\geq 0\ \forall i=1,\ldots,m+1,\
\sum_{i=1}^{m+1} \lambda_i = 1\right\}.$$

Given a function
$\psi\colon\R^m\rightarrow\Rbar$, we shall denote
by $\dom(\psi)$ its effective domain, defined as the subset
of $\R^m$ $\{\xi\suth\psi(\xi)<+\infty\}$,
and by ${\rm epi}\,\psi$
its epigraph, that is the set:
$${\rm epi}\,\psi\dug\{(x,a)\in\R^m\times\R\suth\psi(x)\leq a\}.$$

If $\psi\colon\R^m\rightarrow\Rbar\/$ is Lipschitz continuous in a
neighborhood of a point $\xi$, we shall denote by
$\partial\psi(\xi)$ the generalized gradient
of $\psi$ at $\xi$, defined by
$$\partial\psi(\xi)\dug
{\rm co}\,\left\{
\lim_{i\to +\infty}\nabla\psi(\xi_i)\suth \xi_i\to\xi,\
\xi_i\in\D(\psi)\right\}\,,\nfor{cosb}$$
where $\D(\psi)$ denotes the set of points of differentiability
of $\psi$.
If $\psi$ is differentiable at $\xi$, then
$\partial\psi(\xi)=\{\nabla\psi(\xi)\}$.
We recall that a Lipschitz continuous function $\psi$
is almost everywhere differentiable in $\intdom{\psi}$.

A function $\psi\colon\R^m\rightarrow\Rex$ is convex if, for every
$\xi,\eta\in\R^m$ and for every $\lambda\in[0,1]$, we have
$\psi(\lambda\xi+(1-\lambda)\eta)\leq
\lambda\psi(\xi)+(1-\lambda)\psi(\eta)$.
We say that $\psi$ is concave if $-\psi$
is convex.

Given a function
$\psi\colon\R^m\rightarrow\Rex$, we shall denote by $\psi^{*}$
its dual function, defined for every $p\in\R^m$ by
$$\psi^*(p)\dug\sup_{\xi\in\R^m}\left\{
\scal p{\xi} -\psi(\xi)\right\}.$$
It is well known that the bidual function $\psi^{**}$ coincides
with the convexification of $\psi$, which is
the largest convex function $\varphi$
satisfying $\varphi\leq\psi$.

If $\psi\colon\R^m\to\Rex$ is convex, then the
generalized gradient of $\psi$
coincides in $\intdom{\psi}$
with the
subgradient of $\psi$ in the sense of convex
analysis, defined at every point $\xi\in\dom(\psi)$ by
$$\partial\psi(\xi)\dug\left\{
p\in\R^m\suth \psi(\eta)\geq\psi(\xi)+\scal p{\eta-\xi},\
{\rm for\ every}\ \eta\in\R^m\right\}\nfor{subgr}$$
(see [\rf{Cl}], Proposition~2.2.7).
By definition, we set $\partial\psi(\xi)\dug\emptyset$ for
every $\xi\not\in\dom(\psi)$.

In the following proposition we collect some well
known properties of the subgradient (see [\rf{Cl}] and [\rf{ET}]).

\prop{propsb}{Let $\psi\colon\R^m\to\Rex$ be a convex function.
Then the following properties hold:
\meti{(i)} if $\psi$ is bounded from above in a non--empty
open set $A$, then $\psi$
is locally Lipschitz continuous in $A$;

\meti{(ii)} for every $\xi\in\R^m$,
the set $\partial\psi(\xi)$
(possibly empty) is convex and closed in $\R^m$;

\meti{(iii)} if $\xi\in\intdom{\psi}$, then $\partial\psi(\xi)$
is a non--empty compact set.

%\meti{(iv)} $p\in\partial\psi(\xi)$ if and only if
%$\xi\in\partial\psi^*(p)$;

%\meti{(v)} $p\in\partial\psi(\xi)$ if and only if
%$\psi(\xi)+\psi^*(p)=\scal p{\xi}$.

}

\newpar{clo}{The closure result}

In this section we shall prove a result
concerning the closure of
the convex hull of the epigraph of functions
possibly without superlinear growth.

We recall the notion of strict convexity at infinity,
introduced by Clarke and Loewen in [\rf{CL}].

\dhef{SCI}
{A convex function $\psi\colon\R^m\rightarrow\R$ is
said to be strictly convex at infinity if its graph
contains no rays, that is
for every $\nu\in\R^m$, $\nu\neq 0$,
and for every $\xi\in\R^m$,
the function
$\psx(\vmu)\dug \psi(\vmu\nu+\xi)$ has
the following property:
for every $\vmu_0\in\D(\psx)$
there exists $\vmu_1\in\D(\psx)$,
$\vmu_1>\vmu_0$, such that
$\psx'(\vmu_1)>\psx'(\vmu_0)$.}

\ohss{caratt}
It is easy to see that,
if $\psi\colon\R^m\to\R$ is convex, then
$\psi$ is strictly convex at infinity
if and only if
$\partial\psi^*(p)$ is either empty or
bounded for every $p\in\R^m$.

\dhef{claF}
{We shall denote by $\G$ the family of all
lower semicontinuous functions
$\psi:\R^m\rightarrow\R$ such that
$\psi^{**}\not\equiv -\infty$
and $\psi^{**}$
is strictly convex at infinity.}

\ohss{remuno}
Clearly, every strictly convex function is strictly convex at
infinity. Moreover,
every lower semicontinuous superlinear function
$\psi\colon\R^m\rightarrow\R$ belongs to $\G$.
Indeed, denoting by $\varphi$
the convexification $\psi^{**}$,
for every fixed $\nu$, $\xi\in\R^m$, $\nu\neq 0$,
by (\rf{subgr}) it follows that the inequality
$\scal{\nabla\varphi(\vmu\nu+\xi)}{\vmu\nu}
\geq\varphi(\vmu\nu+\xi)-\varphi(\xi)$
holds for every $s\in\D(\vfx)$.
This implies that
$$\vfx'(\vmu) = \scal{\nabla\varphi(\vmu\nu+\xi)}{\nu}\geq
{\varphi(\vmu\nu+\xi)-\varphi(\xi)\over \vmu},\quad
{\rm for\ every}\ s\in\D(\vfx),\ s>0.$$
Since $\psi$ is superlinear,
the last term tends to $+\infty$ as $\vmu$ goes
to $+\infty$.

\lemma{luno}{For every function $\psi\in\G$, satisfying $\psi\geq 0$
and $\psi(0)=0$,
there exist two positive constants $C$, $\rho$ such that
$\psi(\xi)\geq C |\xi|$ for every $|\xi|>\rho$.}

\proof
We can certainly
assume that $\psi$ is convex,
for if not,
we replace $\psi$ by $\psi^{**}$.
We start by proving that $\psi$ is coercive, that is
$\psi(\xi)\rightarrow +\infty$ as $|\xi|\rightarrow +\infty$.
Since $\psi$ is convex, the sets
$\psi^a\dug\{\xi\in\R^m\suth\psi(\xi)<a\}$ are convex subsets of $\R^m$
for every $a\geq 0$.
By contradiction, suppose that there exists $a>0$ such that $\psi^a$
is unbounded. Since $\psi^a$ is convex, it contains at least one
half line $\{\vmu\nu\suth \vmu\geq 0\}$ for some $\nu\in\R^m$, $\nu\neq 0$.
This means that $\psi_{\nu,0}(\vmu)<a$ for every $\vmu\geq 0$.
Since $\psi_{\nu,0}$ is an absolutely continuous function,
then for every $\tau>0$ we have
$$0\leq \psi_{\nu,0}(\tau)-\psi_{\nu,0}(0)=\int_0^{\tau}
\psi_{\nu,0}'(\sigma)\;d\sigma.$$
Hence, there exists
$\vmu_0\in\D(\psi_{\nu,0})\cap [0,\tau]$
such that $\psi_{\nu,0}'(\vmu_0)\geq 0$.
Since $\psi$ is strictly convex at infinity,
there exists $\vmu_1\in\D(\psi_{\nu,0})$,
$\vmu_1>\vmu_0$, such that
$\psi_{\nu,0}'(\vmu_1)>0$. By the convexity of $\psi_{\nu,0}$
it follows that
$$\psi_{\nu,0}(\vmu)\geq\psi_{\nu,0}(\vmu_1)+
(\vmu-\vmu_1)\psi_{\nu,0}'(\vmu_1),\quad
{\rm for\ every}\ \vmu\geq 0,$$
and this implies that $\lim_{\vmu\to+\infty}\psi_{\nu,0}(\vmu)= +\infty$,
in contradiction with $\psi_{\nu,0}<a$.

Since $\psi$ is coercive, there exist two positive constants
$\rho$, $\delta$ such that:
$$\psi(\eta)\geq\delta,\quad{\rm for\ all\ } |\eta|=\rho.$$
If $|\xi|>\rho$, let us define $\lambda\dug \rho/|\xi|$ and
$\eta\dug \lambda\xi$.
By the convexity of $\psi$, and recalling that $\psi(0)=0$,
we get
$$\psi(\xi)\geq {1\over\lambda} \psi(\eta) =
{\psi(\eta)\over \rho} |\xi|\geq
{\delta\over\rho} |\xi|,$$
so that we have done by choosing $C=\delta/\rho$.
\finedim

We are now in a position to prove the closure result.
The proof
is based on the fact that, if $f$ belongs to the class $\G$,
then for every support hyperplane $r$ of $\fs$,
the function $f-r$ belongs to $\G$.
Applying the estimate of Lemma~\rf{luno} to this function,
we can follow the lines of the proof of Lemma~IX.3.3 in
[\rf{ET}].

\theo{clos}{For every $f\in\G$ the set $\coepif$ is closed.}

\proof
Let $(\xi,a)\in\partial(\coepif)$,
where $\partial S$ denotes the boundary of the set $S$,
and let $r(\eta)\dug \scal{c}{\eta}+d\/$ be
an affine function such that the hyperplane
$H\dug\{(\eta,r(\eta))\}$ weakly separates
$\coepif$ and the point $(\xi,a)$.
Let us define the function
$$\phi(\eta)\dug f(\eta+\xi)-r(\eta+\xi).$$
We have  $\fis(\eta) = \fs(\eta+\xi)-r(\eta+\xi)$, $\fis\geq 0$,
$\fis(0)=0$.
Moreover, for every $\nu\in\R^m$, $\nu\neq 0$,
for every $\eta\in\R^m$ and for every $\vmu\in\D(\fs_{\nu,\xi+\eta})$
we have
$(\fis_{\nu,\eta})'(\vmu)=(\fs_{\nu,\xi+\eta})'(\vmu)-\scal c{\nu}$.
Since $\fs$ is strictly convex at infinity,
then so is $\fis$.
By Lemma~\rf{luno}, there exist two positive constants $C,\rho$ such that
$$\fis(\eta)\geq C|\eta|,\quad {\rm for\ every}\ |\eta|\geq\rho.\nfor{fca}$$

Notice that $(\xi,a)\in\coepif$ if and only if
$(0,0)\in\coepifi$. Moreover, $(\xi,a)\in\partial(\coepif)$
if and only if $(0,0)\in\partial(\coepifi)$.
Hence, to prove the proposition, it suffices to show that
$(0,0)\in\coepifi$.

Let $(\xi^n, a^n)\in\coepifi$ be such that
$\lim_n(\xi^n,a^n)= (0,0)$.
By the characterization (\rf{coA}) of the convex hull,
for every $n$ there exist
$\tilde{\lambda}^n\in E_{m+2}$ and $\xianj\in {\rm epi}\,\phi$,
$j=1,\ldots,m+2$, such that
$$\sumj \lambda^n_j\xianj = (\xi^n, a^n).$$
By the very definition of epigraph it follows that
$$a^n=\sumj\lambda^n_j a^n_j\geq
\sumj\lambda^n_j \phi(\xi^n_j).\nfor{fcaa}$$
Moreover, (\rf{fcaa}) and the fact that $\phi\geq\fis$
imply that
$a^n\geq\sumj\lambda^n_j\fis(\xi^n_j)$.
Since $\fis\geq 0$, the inequality
$$a^n\geq\lambda^n_j\fis(\xi^n_j)\nfor{fcb}$$
holds for every $j=1,\ldots,m+2$.
Let $J\subset\{1,\ldots,m+2\}$ be the
set of all $j$ such that
$\{|\xi^n_j|\}_n$ is unbounded, and let
$I\dug \{1,\ldots,m+2\}\backslash J$.
By passing to a subsequence, we can assume that
there exist $\oxij$, $j\in I$, and $\tilde{\lambda}\in E_{m+2}$,
such that
$$\eqalign{
&\limn|\xi^n_j|=+\infty,\quad j\in J,\cr
&\limn\xi^n_j=\oxij,\quad j\in I,\cr
&\limn\lambda^n_j=\olaj,\quad j\in\{1,\ldots,m+2\}.\cr}$$

For every $j\in J$, we have $|\xi^n_j|>\rho$ for $n$
large enough, and then from (\rf{fca}) and (\rf{fcb}) it follows that
$a^n\geq C\lambda^n_j |\xi^n_j|$.
Since $\lim_n a^n=0$, we get
$$\limn
\lambda^n_j |\xi^n_j|= 0,\quad j\in J.\nfor{aab}$$
{}From (\rf{aab}), and recalling that
$\lim_n\xi^n=0$, we deduce that
$$\eqalign{
\sumI \olaj\oxij&=\limn\sumI \lambda^n_j \xi^n_j =\cr
&=\limn\left(\sumj\lambda^n_j \xi^n_j-
\sumJ\lambda^n_j \xi^n_j\right)=
\limn\left(\xi^n-\sumJ\lambda^n_j \xi^n_j\right)=0.\cr}\nfor{fcc}$$
Moreover, since $\lim_n\lambda^n_j= 0$
for every $j\in J$, we obtain
$$\sumI \olaj=\lim_{n\rightarrow+\infty}
\sumI \lambda^n_j= 1.\nfor{fcd}$$
Since $\phi$ is a non--negative
lower semicontinuous function, we get
$$0\leq\sumI \olaj\phi(\oxij)\leq
\liminf_{n\rightarrow +\infty}\sumI \lambda^n_j\phi(\xi^n_j)
\leq\liminf_{n\rightarrow+\infty} a^n =0.\nfor{aac}$$
There is no loss of generality in assuming
that $\lambda_j>0$
for every $j\in I$, hence (\rf{aac}) implies
that $\phi(\oxij)=0$ for every $j\in I$,
that is $(\oxij,0)\in{\rm epi\,}\phi$ for every $j\in I$.
Thus, by (\rf{fcc}) and (\rf{fcd}),
we can conclude that
$(0,0)$ belongs to $\coepifi$.
\finedim

Now we state two direct consequences of Theorem~\rf{clos}.

\corol{coruno}{If $f\in\G$, then
$$\fs(\xi)=\min\left\{\sum_{j=1}^{m+1} \lambda_j f(\xi_j)\suth
\sum_{j=1}^{m+1}\lambda_j\xi_j=\xi,\
\tilde{\lambda}\in E_{m+1}\right\},$$
for every $\xi\in\R^m$.}

\proof
See [\rf{ET}], Lemma~IX.3.3.
\finedim

We recall that a function $f\colon I\times\R^m\to\R$
is said to be a normal integrand (see [\rf{ET}])
if $f(t,\cdot)$ is lower semicontinuous for a.e.~$t\in I$,
and there exists a Borel function
$\tilde{f}\colon I\times\R^m\to\R$ such that
$\tilde{f}(t,\cdot)=f(t,\cdot)$ for a.e.~$t\in I$.

\corol{cordue}{
Let $f\colon I\times\R^m\to\R$  be a normal integrand, and
suppose that $f(t,\cdot)\in\G$ for every $t\in I$.
Then for any measurable mapping
$p\colon[0,T]\rightarrow\R^m$, there exist a measurable mapping
$\tilde{\lambda}\colon[0,T]\rightarrow E_{m+1}$ and
$m+1$ measurable mappings $q_j\colon[0,T]\rightarrow\R^m$,
such that
$$\sum_{j=1}^{m+1}\lambda_j(t) q_j(t)=p(t),\quad
\sum_{j=1}^{m+1}\lambda_j(t) f(t, q_j(t)) = \fs(t, p(t)),$$
for almost all $t\in [0,T]$.}

\proof
See [\rf{ET}], Proposition~IX.3.1.
\finedim
\newpar{exi}{Existence results for variational problems}

In this section we shall show that the existence result
proved by Cellina and Colombo in [\rf{CC}]
holds even for functions of the class $\E$ defined below.
In the following, the convexification and the gradient
of a function $\psi(t,\xi)$
are understood with respect to $\xi$.

\dhef{claE}
{We shall denote by $\E$ the family of all functions
$\psi\colon I\times\R^m\rightarrow\R$, bounded from below, such that
$\psi(\cdot,\xi)$ is Lipschitz continuous for every
fixed $\xi\in\R^m$,
$\psi(t,\cdot)$ is lower semicontinuous for every
fixed $t\in I$, and
$$\lim_{R\to+\infty}
\sup_{\scriptstyle t\in I\atop\scriptstyle |\xi|>R}
\sup\left\{\psi^{**}(t,\xi)-\scal p{\xi}\suth
p\in\dexi\psi^{**}(t,\xi)\right\}=-\infty.\nfor{cen}$$}

The following proposition gives a characterization
of the family $\E$. The proof is similar to the one
of Proposition~3.2 in [\rf{Cra}].

\prop{carE}{The condition (\rf{cen}) in Definition~\rf{claE}
is equivalent to:
$$\limn\left[\psin-\scal{\nabla\psin}{\xi^n}\right]=-\infty\nfor{cend}$$
for every sequence $(t^n,\xi^n)\in I\times\R^m$ such that
$\xi^n\in\D(\psi^{**}(t^n,\cdot))$, $\lim_n|\xi^n|=+\infty$.}

\proof
We have to prove that (\rf{cend}) implies (\rf{cen}),
the other implication being trivial.
Let us denote by $\chi(R)$ the argument of the limit in (\rf{cen}),
and let $\{R_n\}$ be a diverging sequence.
For every fixed $n\in\nat$,
by definition of supremum, there exists
$(t^n,\xi^n,p^n)\in I\times\R^m\times\R^m$, with
$p^n\in\dexi\psin$ and $|\xi^n|>R_n$, such that
$$\chi(R_n)\leq\psin-\scal{p^n}{\xi^n}+1.\nfor{sta}$$
{}From (\rf{cosb}) and
(\rf{coA}), there exist
$p^n_j\in\dexi\psin$,
$\xi^n_j\in\D(\psi^{**}(t^n,\cdot))$,
with $|\xi^n_j-\xi^n|<1$,
$j\in J\dug\{1,\ldots,m+1\}$,
and $\tilde{\lambda}^n\in E_{m+1}$, such that
$$p^n=\sj\lnj p^n_j,\qquad
\left|\nabla\psinj - p^n_j\right|<{1\over |\xi^n|+1},\quad
{\rm for\ every}\ j\in J.$$
For every $j\in J$, the last inequality
and the fact that $|\xi^n_j-\xi^n|<1$
imply that
$$\left|\scal{\nabla\psinj-p^n_j}{\xi^n_j}\right|<
{|\xi^n_j|\over |\xi^n|+1}<1\,.\nfor{stb}$$
By the convexity of $\psi^{**}(t^n,\cdot)$ we have
$$\psin-\psinj\leq\scal{p^n_j}{\xi^n-\xi^n_j},\quad
{\rm for\ every}\ j\in J.\nfor{stc}$$
Using (\rf{stb}) and (\rf{stc}) we obtain
$$\psin-\scal{p^n_j}{\xi^n}\leq
\psinj-\scal{\nabla\psinj}{\xi^n_j}+1\,.\nfor{std}$$
Multiplying (\rf{std}) by $\lnj$ and summing over $j$
it follows that
$\psin-\scal{p^n}{\xi^n}\leq\mu^n$,
where $\mu^n\dug1+\max_j\{\psinj-\scal{\nabla\psinj}{\xi^n_j}\}$.

Since $\lim_n|\xi^n_j|=+\infty$ for every $j\in J$,
(\rf{cend}) implies that $\lim_n\mu^n=-\infty$.
Hence, by (\rf{sta}), it follows that
$$\limn \chi(R_n)\leq\limn (\mu^n+1) =-\infty\,.$$
Since $\chi$ is a monotone non--increasing function,
(\rf{cen}) holds.
\finedim

\ohss{remzero}
The Definition~\rf{claE}
agrees with the one given in [\rf{CTZ}] and [\rf{Cra}],
respectively
in the case of
convex time--independent smooth functions and
non--convex time--independent functions.

\lemma{EG}{If $\psi\in\E$, then
$\psi(t,\cdot)\in\G$ for every $t\in I$.}

\proof
Let us fix $t\in I$, and denote by $\varphi$ the convexification
with respect to $\xi$ of $\psi(t,\xi)$.
By Lemma~3.3 in [\rf{Cra}], the effective domain $\dom(\varphi^*)$
of $\varphi^*$ is an open subset of $\R^m$.
Hence, by Proposition~\rf{propsb}(iii),
$\partial\varphi^*(p)$ is either bounded, if $p\in\dom(\varphi^*)$,
or empty, if $p\not\in\dom(\varphi^*)$.
By Remark~\rf{caratt}, the result is thus proved.
\finedim

\lemma{esti}{Let $\varphi\colon I\times\R^m\times\R^m\rightarrow\R$
be a lower semicontinuous function,
Lipschitz continuous with respect to the first variable.
Assume that $\varphi(t,x,\cdot)$ is convex for a.e.~$t\in I$
and for every $x\in\R^m$, and that there exist three constants
$C_i$, $i=0,1,2$, such that
$$
|v|\le C_0|\varphi(t,x,\xi)|+C_1|x|+C_2\,,
\nfor{hfitre}$$
for every
$(t,x,\xi)\in I\times\R^m\times\R^m$ and for every
$v\in\partial_t\varphi(t,x,\xi)$,
where $\partial_t\varphi$ denotes the generalized gradient
of $\varphi$ with respect to $t$.

Let $u\in W^{1,1}(I,\R^m)$, and
assume that the function $t\mapsto \varphi(t,u(t),u'(t))$ belongs
to $L^1(I)$.
Then there exists $k_0\in L^1(I)$ such that
$$\left|\varphi(s_2,u(t),u'(t))-\varphi(s_1,u(t),u'(t))\right|
\leq k_0(t) |s_2-s_1|,$$
for every $t,s_1,s_2\in I$.}

\proof
For every fixed $t_1$, $t_2\in I$, let us define the function
$$g(s)\dug\left|\varphi(t_1+s d,x,\xi)-\varphi(t_1,x,\xi)\right|,
\quad s\in [0,1],$$
where $d\dug t_2-t_1$.
By (\rf{hfitre}), it follows that for a.e.~$s\in[0,1]$
$$g'(s)\leq |d||\ptf(t_1+s d,x,\xi)|\leq
|d|\left(C_0 g(s) +C_0 |\varphi(t_1,x,\xi)| + C_1 |x|+C_2\right).$$
We can apply Gronwall's inequality
to the non--negative absolutely continuous function $g$,
obtaining
$$|\varphi(t_2,x,\xi)-\varphi(t_1,x,\xi)|=g(1)\leq
|t_2-t_1|\ect\left(C_0 |\varphi(t_1,x,\xi)|+C_1 |x| + C_2\right).
\nfor{beuno}$$
This inequality, with $t_1=t$ and $t_2=s_1$, implies that
$$|\varphi(s_1,x,\xi)|\leq |\varphi(t,x,\xi)|+
T\ect\left(C_0 |\varphi(t,x,\xi)|+C_1|x|+C_2\right).
\nfor{bedue}$$
Again by (\rf{beuno}), with $t_1=s_1$, $t_2=s_2$, and by
(\rf{bedue}), it follows that
$$
|\varphi(s_2,x,\xi)-\varphi(s_1,x,\xi)|\leq
|s_2-s_1| (\tilde{C}_0|\varphi(t,x,\xi)|+
\tilde{C}_1|x|+\tilde{C}_2),$$
where $\tilde{C}_i\dug C_i\ect (1+T C_0\ect)$, $i=0,1,2$.
Finally, by hypothesis, the function
$$k_0(t)\dug \tilde{C}_0|\varphi(t,u(t),u'(t))|+
\tilde{C}_1|u(t)|+\tilde{C}_2$$
belongs to $L^1(I)$, completing the proof.
\finedim

\dhef{nag}{We shall say that
$\theta\in\C^1((0,+\infty),\R)$
is a Nagumo function if $\theta$
is convex, increasing and it satisfies
$\lim_{r\rightarrow +\infty}{{\theta(r)}/r}=+\infty$.}

We begin the study of minimization problems,
starting with
an existence result for convex functionals.
We collect here the basic hypotheses on the integrand.

\meti{$(H_0)$} $f\in\E$, and $f(t,\cdot)$ is a convex function
for every $t\in I$.

\meti{$(H_1)$} There exist two constants $A$ and $B$,
with $B>0$, such that
$f(t,\xi) \ge -A+B|\xi|$
for every $(t,\xi)\in I\times\R^m$.

\meti{$(H_2)$} $g\colon I\times\R^m\rightarrow\R$
is Lipschitz continuous with respect to the first variable,
continuous with respect to the second,
and there exist two constants $\alpha$, $\beta$,
with $0\le\beta<B/T$, such that $g(t,x)\ge -\alpha-\beta|x|$
for every $(t,x)\in I\times\R^m$.

\meti{$(H_3)$} There exist three constants $C_i$,
$i=0,1,2$, such that the condition
(\rf{hfitre}) holds with
$\varphi(t,x,\xi)\dug g(t,x)+f(t,\xi)$.

\ohss{bound}
If $f\in\E$ is independent of $t$, then
it is easily seen that
Lemma~\rf{luno} and Lemma~\rf{EG}
imply that condition $(H_1)$ is always satisfied
for suitable constants $A$, $B$, with $B>0$.

\def\solu{\tilde{u}}
\theo{app}{
Let $f$ and $g$ satisfy the hypotheses $(H_0)$, $(H_1)$, $(H_2)$,
$(H_3)$.
Then there exists a solution $\solu$ to the problem
$$
\min\left\{F(u)\ \big|\
u\in W^{1,1}(I,\R^m),\, u(0)=a,\, u(T)=b\right\}\nfor{min}
$$
where
$$
F(u)\dug\int_I [f(t,u'(t))+g(t,u(t))]\, dt\,.
$$
Moreover $\solu$ belongs to $W^{1,\infty}(I, \R^m)$ and satisfies
for a.e.~$t\in I$
$$
f(t,\solu'(t))-\scal{p(t)}{\solu'(t)}
+g(t,\solu(t))=c +\int_0^t v(\tau)\, d\tau\,,\nfor{DR}
$$
where $c$ is a constant, and
$(v(t),p(t))\in
(\partial_t f(t,\solu'(t))+\partial_t g(t,\solu(t)),
\dexi f(t,\solu'(t)))$
for almost every $t\in I$.}

\proof The proof follows the lines of the one of Theorem~3 in [\rf{Cla}],
with some changes due to the fact that in this case the lagrangian
is not bounded from below. As in [\rf{Cla}] one can prove, using
the De Giorgi's semicontinuity result (see [\rf{Butt}])
and the Dunford--Pettis criterion
of weak compactness in $L^1(I, \R^m)$, that for every Nagumo function
$\theta$ and for every $l>0$ there exists a solution $u_l$ to the
problem
$$
\min\left\{F(u)\ \big|\
u\in AC_\theta^l(I,\R^m),\, u(0)=a,\, u(T)=b\right\}\,,
$$
where $AC_\theta^l(I,\R^m)$ denotes the class of all function $u \in
W^{1,1}(I,\R^m)$ such that $\Theta(u)\le l$, with
$\Theta(u)\dug\int_I \theta(|u'(t)|)\, dt\,.$
Let us set $V_\theta(l)\dug F(u_l)$.

One can easily check that, if $V_\theta(l)=V_\theta(l_0)$ for every
$l\ge l_0$, then $u_{l_0}$ is a solution to the problem
$$
\min\left\{F(u)\ \big|\
u\in W^{1,1}(I,\R^m),\, \Theta(u)<+\infty,\ u(0)=a,
\, u(T)=b\right\}\,.\nfor{pteta}
$$
Finally, as in [\rf{Cla}], if we are able to prove that $u_{l_0}$ belongs
to $W^{1,\infty}(I,\R^m)$, then we can conclude that such a function
is a solution to (\rf{min}).

Thus it remains to prove that $V_\theta$
is eventually constant and that, for $l$ large enough, $u_l$ belongs
to $W^{1,\infty}(I,\R^m)$ and satisfies (\rf{DR}).
Since $V_\theta$ is lower semicontinuous,
for every $l>0$ there exists a proximal subgradient
(see [\rf{Clar}]) of $V_\theta$ at
$l$ and, since $V_\theta$ is nonincreasing,
it is nonpositive.
If $V_\theta$ is not eventually constant,
by Proposition~6.1 in [\rf{Cla}],
there exists a diverging sequence $\{l_k\}$ such that
the proximal subgradient of $V_\theta$
at $l_k$ takes the form $-r_k$, with $r_k>0$.
Moreover, it is easy to check that,
if we set $u_k\equiv u_{l_k}$, then
$\Theta(u_k)=l_k$, so that
$$
\lim_{k \rightarrow +\infty}
\|u_k'\|_{L^{\infty}}\ge \lim_{k \rightarrow +\infty}
\theta^{-1}(l_k/T)=+\infty\,.\nfor{unb}
$$
By definition of $r_k$ and the fact that $\Theta(u_k)=l_k$,
it follows that
for every $k\in\nat$ there exists
a positive constant $\sigma_k$ such that, if we define
$$
G(u)\dug F(u)+r_k \Theta(u) +\sigma_k |\Theta(u)-\Theta(u_k)|^2\,,
$$
then we get that $G(u_k)\le G(u)$ for every
$u$ admissible for (\rf{pteta})
and such that $\Theta(u)$ is sufficiently near to $\Theta(u_k)$
(see [\rf{Cla}]).
By $(H_3)$ and Lemma~\rf{esti}, it follows that there exists $k_0
\in L^1(I)$ such that for every $s_1, s_2, t\in I$
$$
|f(s_1,u_k'(t))+g(s_1,u_k(t))-f(s_2,u_k'(t))-g(s_2,u_k(t))|
\le k_0(t)|s_1-s_2|\,,
$$
so that we can apply Theorem 5 of [\rf{Cla}]. Thus we obtain that
$u_k$ satisfies
$$
E_f(t,u_k'(t))+g(t,u_k(t))+r_k E_\theta(|u_k'(t)|)
=c_k +\int_0^t v_k(\tau)\, d\tau\,,\nfor{DRk}
$$
where $c_k$ is a constant,
$E_f(t,u_k'(t))\dug f(t,u_k'(t))-\scal{p_k(t)}{u_k'(t)}$,
with $(v_k(t),p_k(t))\in (\partial_t f(t,u_k'(t))
+\partial_t g(t,u_k(t)),\dexi f(t,u_k'(t)))$
for a.e.~$t\in I$, and
$E_{\theta}(\vmu)\dug \theta(\vmu)-\vmu \theta'(\vmu)$.

Moreover there exists $M_1>0$ such that $\|u_k\|_{L^\infty}
\le M_1$ for every $k\in\nat$.
Actually, if there exists $t_k\in I$ such that
$\limsup_k |u_k(t_k)|=+\infty$, then
$$
\limsup_{k\rightarrow +\infty}\int_I |u_k'(t)|\, dt \ge
\limsup_{k\rightarrow +\infty} \left|\int_0^{t_k} u_k'(t)\, dt\right|=
\limsup_{k\rightarrow +\infty}|u_k(t_k)-a|=+\infty\,,
$$
while, if we define $u_0(t)\dug a+\xi t$, with $\xi\dug (b-a)/T$,
then $u_0$ is admissible for (\rf{pteta}), $F(u_0)<+\infty$, and
$$
F(u_0)\ge F(u_k) \ge (-A-\alpha)T+B\|u_k'\|_{L^1}-\beta \|u_k\|_{L^1}
\ge \tilde{A}+(B-\beta T)\|u_k'\|_{L^1}\,,\nfor{bdaa}
$$
so that, by $(H_2)$, $\{u_k'\}$ must be bounded in $L^1(I,\R^m)$.

The boundedness
of $\{u_k\}$ in $L^\infty(I,\R^m)$ and the continuity of $g$
guarantee that there exists $M_2$ such that
$$
|g(t,u_k(t)| \le M_2\,,\nfor{cip}
$$
for a.e.~$t\in I$ and for every $k$. Moreover, by $(H_3)$ we obtain
$$
\eqalign{&\left|\int_0^t  v_k(s)\, ds\right| \le
\int_I\left[C_0 \left|f(s,u_k'(s))+g(s,u_k(s))\right|
+C_1|u_k(s)| +C_2\right]\, ds\leq\cr
&\le \int_I \left[C_0 \left|
\alpha+\beta|u_k(s)|+f(s,u_k'(s))+g(s,u_k(s))
\right|+\tilde{C}_1|u_k(s)| +\tilde{C}_2
\right]\, ds\,, \cr}\nfor{miao}
$$
where $\tilde{C}_1\dug C_0\beta+C_1$ and
$\tilde{C}_2\dug C_0 |\alpha|+C_2$.
Without loss of generality we can assume that $f$ is positive, so
that, thanks to $(H_2)$, it follows that for every $k\in\nat$
$$
f(s,u_k'(s))+g(s,u_k(s))+\alpha+
\beta|u_k(s)| \ge 0,\quad {\rm a.e.}\ s\in I\,.\nfor{posi}
$$
By (\rf{bdaa}), (\rf{miao}) and (\rf{posi}) there exist
$M_3>0$ and two constants $\hat{C}_1$, $\hat{C}_2$ such that
$$
\left|\int_0^t v_k(s)\, ds\right| \le
C_0 F(u_k)+\hat{C}_1\|u_k\|_{L^1} + \hat{C}_2\le M_3\,,
\quad {\rm for\ every}\ t\in I. \nfor{bau}
$$
By (\rf{DRk}), (\rf{cip}), and (\rf{bau}) we obtain
$$
E_f(t,u_k'(t))+r_kE_\theta(|u_k'(t)|)
\le c_k +M_2+M_3\,,
$$
for every $t\in I$ and for every $k\in\nat$.

We claim that it is not possible that there exists a
subsequence of $\{c_k\}$, still denoted by $\{c_k\}$,
such that
$\lim_k c_k=-\infty$.
Indeed, if this is the case,
then for every $t\in I$ we should have
$$
\lim_{k\rightarrow +\infty}
E_f(t,u_k'(t))+r_kE_\theta(|u_k'(t)|)=-\infty\,.\nfor{enin}
$$
Since $f\in\E$ and $\theta$ is superlinear, (\rf{enin}) implies that
$\lim_k|u_k'(t)|=+\infty$ for every $t\in I$, which, by Fatou's
Lemma, contradicts the
boundedness of $u_k'$ in $L^1(I,\R^m)$.

Thus there exists $c^*$ such
that $c_k \ge c^*$ for every $k$. From (\rf{DRk}) we obtain,
for every $t\in I$,
$$
E_f(t,u_k'(t))+r_k E_\theta(|u_k'(t)|)\ge c^*-M_2-M_3\,.\nfor{ben}
$$
Now let us suppose that for every $k$ there exists $t_k\in I$ such that
$\limsup_k |\xi_k|=+\infty$, where $\xi_k\dug u_k'(t_k)$.
Since $f$ and $\theta$ belong to $\E$,
we have
$$
\liminf_{k\to+\infty} \left[E_f(t_k,\xi_k)+
r_kE_\theta(|\xi_k|)\right]\leq
\liminf_{k\to+\infty} \sup_{t\in I} \left\{
E_f(t,\xi_k)+
r_kE_\theta(|\xi_k|)\right\} = -\infty\,,
$$
in contradiction with (\rf{ben}).
This implies that $\|u_k'\|_{L^\infty}$
is bounded, which contradicts (\rf{unb}).

So we can conclude that $V_\theta$ is eventually constant. Hence
for $k$ sufficiently large $u_k\in W^{1,\infty}(I, \R^m)
$ is a solution of (\rf{pteta}).
Moreover $r_k=0$, so that $u_k$ satisfies (\rf{DR}). Then the proof
is complete. \finedim

The last part of this section is devoted to the study of the
non--convex case.
The hypotheses $(H_0)$ and $(H_3)$ will be replaced respectively by:

\meti{$(H_0')$} $f\in\E$.

\meti{$(H_3')$} There exist three constants $C_i$,
$i=0,1,2$, such that the condition
(\rf{hfitre}) holds with
$\varphi(t,x,\xi)\dug g(t,x)+\fs(t,\xi)$.

Notice that $(H_3')$ requires the Lipschitz continuity of
$\fs$ with respect to $t$.
The following two lemmas show that this conclusion
follows from $(H_0')$ and

\meti{$(H_4)$} For every $R>0$ there exists a constant $L$
such that
$$|f(t,\xi)-f(s,\xi)|\leq L|t-s|,\quad
{\rm for\ every}\ t,s\in I,\ {\rm and}\ \xi\in \overline{B}_R,$$

where $\overline{B}_R$ denotes the closed ball
centered at the origin and with radius $R$.

\lemma{cca}
{Let $\psi\in\E$, and let us define, for every $(t,p)\in I\times\R^m$,
the set
$$W(t,p)\dug \{\xi\in\R^m\suth p\in\dexi\psi^{**}(t,\xi)\}\,.$$
Then for every $r>0$ there exists $R>0$ such that
for every $(t,p)\in I\times\R^m$
the condition $W(t,p)\cap \overline{B}_r\neq\emptyset$
implies $W(t,p)\subset\overline{B}_R$.}

\proof
Suppose, by contradiction, that there exist sequences
$(t_n, p_n)\subset I\times\R^m$, $(\eta_n)\subset\overline{B}_r$,
$(\xi_n)\subset\R^m$, with
$\lim_n |\xi_n|=+\infty$, such that, for every $n\in\nat$,
$$p_n\in\dexi\psi^{**}(t_n,\eta_n),\quad
p_n\in\dexi\psi^{**}(t_n,\xi_n)\,.\nfor{ccfa}$$
{}From (\rf{ccfa}) it follows that, for every $n\in\nat$,
$$\psi^{**}(t_n,\eta_n)-\scal{p_n}{\eta_n} =
\psi^{**}(t_n,\xi_n)-\scal{p_n}{\xi_n}\,.\nfor{ccfb}$$
Since $(\eta_n)$ is a bounded sequence, there exists a constant
$C$ such that the left hand side of (\rf{ccfb}) is
bounded from below by $C$. Thus
$$C\leq \psi^{**}(t_n,\xi_n)-\scal{p_n}{\xi_n}
\leq\chi(|\xi_n|),\quad {\rm for\ every}\ n\in\nat,\nfor{ccfc}$$
where $\chi(R)$ is the argument of the limit in
(\rf{cen}).
Since $\lim_n |\xi_n|=+\infty$, from (\rf{cen}) we have
that $\lim_n \chi(|\xi_n|)=-\infty$, which contradicts (\rf{ccfc}).
\finedim

\ohss{ccr}
{Let us fix $\xi\in\R^m$.
Let $t\in I$, $\tilde{\lambda}\in E_{m+1}$, $\xi_j\in\R^m$,
$j=1,\ldots,m+1$ satisfy
$$\fs(t,\xi)=\sum_{j=1}^{m+1}\lambda_j f(t,\xi_j),\quad
\xi=\sum_{j=1}^{m+1}\lambda_j \xi_j.$$
Since for every $j$ there exists $p_j\in\dexi\fs(t,\xi)$
such that $\xi_j\in W(t,p_j)$, by Lemma~\rf{cca} we
obtain that there exists $R>0$, depending only on $|\xi|$,
such that $\xi_j\in\overline{B}_R$ for every
$j=1,\ldots,m+1$.}

\lemma{ccb}
{If $f\in\E$ satisfies $(H_4)$, then $\fs(\cdot, \xi)$ is
Lipschitz continuous for every $\xi\in\R^m$.}

\proof
Let us fix $\xi\in\R^m$, and consider $t$, $s\in I$.
By Corollary~\rf{coruno}, there exist
$\tilde{\lambda}$, $\tilde{\mu}\in E_{m+1}$,
$\xi_j$, $\eta_j\in\R^m$, $j=1,\ldots,m+1$,
such that
$$\fs(t,\xi)=\sum_{j=1}^{m+1}\lambda_j f(t,\xi_j),\quad
\fs(s,\xi) = \sum_{j=1}^{m+1}\mu_j f(s,\eta_j)\,,$$
and $\xi=\sum_j\lambda_j\xi_j=\sum_j\mu_j\eta_j$.
Moreover, one has
$$\fs(t,\xi)\leq\sum_{j=1}^{m+1}\mu_j f(t,\eta_j),\quad
\fs(s,\xi) \leq \sum_{j=1}^{m+1}\lambda_j f(s,\xi_j)\,.$$
Then, by Remark~\rf{ccr} and $(H_4)$, there exists $L>0$,
depending only on $|\xi|$, such that
$$\fs(s,\xi)-\fs(t,\xi)\leq
\sum_{j=1}^{m+1}\lambda_j [f(s,\xi_j)-f(t,\xi_j)]\leq
\sum_{j=1}^{m+1}\lambda_j L|t-s|= L|t-s|\,.$$
In the same way one obtains
$$\fs(t,\xi)-\fs(s,\xi)\leq
\sum_{j=1}^{m+1}\mu_j [f(t,\eta_j)-f(s,\eta_j)]\leq L|t-s|\,,$$
completing the proof.
\finedim

We are now in a position to prove the existence result for the
non--convex case.

\theo{conca}
{Let $g$ and $f$ satisfy the basic hypotheses
$(H_0')$, $(H_1)$, $(H_2)$, $(H_3')$, $(H_4)$,
and assume that
$g(t,\cdot)$ is concave for every $t\in I$.
Then the problem (\rf{min}) has a solution
$u\in W^{1,\infty}([0,T],\R^m)$.}

\proof
The proof follows the same lines of the one
of Theorem~1 in [\rf{CC}].
It is enough to use Theorem~\rf{app}
to obtain a solution
$\tilde{u}\in W^{1,\infty}([0,T],\R^m)$
of the relaxed problem,
and to replace
Lemma~IX.3.3 and Proposition~IX.3.1 of [\rf{ET}]
with Corollaries~\rf{coruno} and \rf{cordue}.
Since $\tilde{u}'\in L^{\infty}([0,T],\R^m)$, it is easily
seen that we obtain a solution
$u\in W^{1,\infty}([0,T],\R^m)$.
\finedim

%
%
%%%%%%%
%
%
% bibliografia
%
\bigskip\medskip\noindent
\vfill\eject
{\bf References}\medskip\goodbreak

\bibart{AAB}{Ambrosio, L., Ascenzi, O.~and Buttazzo, G.}
{Lipschitz regularity for minimizers of integral functionals
with highly discontinuous integrands}
{J.~Math.~Anal.~Appl.}{142}{1989}{301}{316}

%\biblib{AC}{Aubin, J.P.~and Cellina, A.}
%{Differential Inclusions}
%{Springer--Verlag}{Berlin}{1984}

\bibart{BD}{Botteron, B.~and Dacorogna, B.}
{Existence and non-existence results for non-coercive
variational problems and applications in ecology}
{J.~Differential Equations}{85}{1990}{214}{235}

\bibart{BM}{Botteron, B.~and Marcellini, P.}
{A general approach to the existence of minimizers of one
dimensional non-coercive integrals of the calculus of variations}
{Ann.~Inst.~H. Poin\-ca\-r\'e}{8}{1991}{197}{223}

\biblib{Butt}{Buttazzo, G.}
{Semicontinuity, Relaxation and Integral Representation
in the Calculus of Variations}
{Pitman Res.~Notes Math.~Ser.~{\bf 207}}
{Longman, Harlow}{1989}

\bibart{CC}{Cellina, A.~and Colombo, G.}
{On a classical problem of the calculus of variations
without convexity assumptions}
{Ann.~Inst.~H.~Poincar\'e}{7}{1990}{97}{106}

\bibprep{CTZ}{Cellina, A., Treu, G.~and Zagatti, S.}
{On the minimum problem for a class of non-coercive functionals}
{Preprint SISSA}{1994}

\biblib{ces}{Cesari, L.}
{Optimization -- theory and applications}
{Springer--Verlag}{New York}{1983}

\biblib{Cl}{Clarke, F.~H.}
{Optimization and Nonsmooth Analysis}
{Wiley Interscience}{New York}{1983}

\biblib{Clar}{Clarke, F.~H.}
{Methods of dynamic and nonsmooth optimization}
{CBMS-NSF Regional Conf.~Ser.~in Appl.~Math.,
Vol.~57, SIAM}{Philadelphia}{1989}

\bibart{Cla}{Clarke, F.~H.}
{An indirect method in the calculus of variations}
{Trans.~Amer.~Math. Soc.}{336}{1993}{655}{673}

\bibart{CL}{Clarke, F.~H.~and Loewen, P.~D.}
{An intermediate existence theory in the Calculus of Variations}
{Ann.~Scuola Norm.~Sup.~Pisa Cl.~Sci.~(4)}{16}{1989}{487}{526}

\bibprep{Cra}{Crasta, G.}
{An existence result for non-coercive non-convex problems
in the Calculus of Variations}
{Preprint SISSA}{1994}

%\biblib{Dei}{Deimling, K.}
%{Multivalued Differential Equations}
%{De Gruyter}{Berlin}{1992}

\biblib{ET}{Ekeland, I.~and Temam, R.}
{Convex Analysis and Variational Problems}
{North--Holland}{Amsterdam}{1977}

\bibart{Mar}{Marcellini, P.}
{Alcune osservazioni sull'esistenza del minimo di integrali
del calcolo delle variazioni senza ipotesi di convessit\`a}
{Rend. Mat.}{13}{1980}{271}{281}

\biblib{Ol}{Olech, C.}
{The Lyapunov theorem: its extensions and applications}
{in ``Methods of Non-convex Analysis'', A.~Cellina ed.}
{Springer--Verlag, New York}{1990}

\bibart{Ray}{Raymond, J.P.}
{Existence theorems in optimal control problems without
convexity assumptions}
{J.~Optim.~Theory Appl.}{67}{1990}{109}{132}

\prebye
\bye